\documentclass[12pt]{amsart}

\usepackage{graphicx}
\usepackage{amsmath}

\usepackage{tikz-cd}

\def\frak{\mathfrak}

\def\C{\mathbb{C}}

\def\P{\mathbb{P}}
\def\R{\mathbb{R}}

\def\Z{\mathbb{Z}}

\def\al{\alpha}

\def\ka{\kappa}

\def\om{\omega}
\def\Ga{\Gamma}

\def\Up{\Upsilon}

\newcommand{\der}{{\rm d}}

\numberwithin{equation}{section}

\newtheorem{theorem}{Theorem}[section]

\newtheorem{proposition}[theorem]{Proposition}

\theoremstyle{remark}

\theoremstyle{remark}

\usepackage{amssymb,stmaryrd}
\usepackage{amscd}

\makeatletter
\@namedef{subjclassname@2020}{%
  \textup{2020} Mathematics Subject Classification}
\makeatother

\author{Matthew Randall}
\address{School of Mathematics and Statistics\\
Nanjing University of Information Science and Technology\\
219 Ningliu Road\\
Nanjing, China}
\email{100093@nuist.edu.cn}

\title{Contact structures of type $G_2$ associated to solutions of Noth's equation}


\begin{document}

\begin{abstract}
We establish a correspondence between solutions of Noth's equation, a non-linear ordinary differential equation that shows up in the theory of $(2,3,5)$-distributions, and diffeomorphisms of any contact structure of type $G_{2}$ to the standard one. 
\end{abstract}

\maketitle

\pagestyle{myheadings}
\markboth{Randall}{Contact structures of type $G_2$ associated to solutions of Noth's equation}

\section{Introduction}

The theory of Lie contact structures of type $G_2$ has origins in Cartan's 1910 paper \cite{Cartan1910} and was more recently investigated in \cite{ky06}, \cite{LNS}, \cite{ufo1}. Let $G_2$ be the split real form of the smallest exceptional simple Lie group $G_2$ and $\frak{g}_2$ be its Lie algebra. Let $P_2$ be the 9-dimensional parabolic subgroup associated to the contact gradation of $\frak{g}_2$ of depth 2. In Proposition 4.4 of \cite{ufo1}, the $\frak{g}_2$ Lie algebra of vector field is presented as a symmetry algebra of a Lie contact structure on a 5-dimensional contact manifold $G_2/P_2$.

Also, a sixth-order ordinary differential equation called Noth's equation, given by 
\begin{align}\label{noth}
10H''^3H^{(6)}-70H''^2H'''H^{(5)}&-49H''^2H''''^2+280H''H'''^2H''''-175H'''^4=0,
\end{align}
shows up as a condition for maximal symmetry for certain types of $(2,3,5)$-distributions \cite{annur}. The solutions to this equation up to $SL(2,K)$ equivalence, where $K=\R$ or $\C$, were found in \cite{r16} and \cite{r17} and can be given by solutions of the generalized Chazy equation with parameter $k=\frac{3}{2}$. The 5-dimensional contact manifold is related to $(2,3,5)$-distributions through the double fibration. We review later on in this manuscript the transformation of this double fibration as given explicitly in Cartan's 1910 paper, and we complete the computation initiated in page 180 of \cite{Cartan1910}, where we give the Lie contact structure of type $G_2$ and the $\frak{g}_2$  Lie algebra of vector fields associated to solutions of Noth's equation, or equivalently the generalized Chazy equation with parameter $\frac{3}{2}$.

The solution $H(t)=3t^2$ to Noth's equation can be seen as the simplest one, and can be related to the standard contact structure of  Proposition 4.4 of \cite{ufo1} through the double fibration. In this paper we give analogous results to Proposition 4.4 of \cite{ufo1}, replacing the standard contact structure corresponding to the solution $H(t)=3t^2$ by the more general solution obtained in \cite{r17}. This gives contact manifolds that are related to the standard contact one through a coordinate diffeomorphism, which induces the corresponding isomorphism of the linear differential system associated to the contact structure, denoted by $\mathcal{C}$, and the symmetric third power of a rank 2 bundle $V$ on the standard contact manifold. We also need to fix the fibre parameter to identify it with the independent variable in Noth's equation, which leads to the following main result of the paper. 

\begin{theorem}\label{maintheorem}
There is a one-to-one correspondence between $SL(2,K)$-equivalent solutions of Noth's equation and isomorphisms of $S^3(V)$ with $\mathcal{C}$ with specified fibre parameter.
\end{theorem}
 
What is striking about the main result is that on one side of the correspondence, we have some solution of an integrable non-linear differential equation of high order and on the other side, we have a completely geometric characterization involving diffeomorphisms of contact structures of type $G_2$. The reader should also note that the solutions of Noth's equation here do not see the Schwarz triangle functions solutions obtained in \cite{r17}, but only the contact symmetry. We first review the background material on contact structures of type $G_2$ in Section \ref{background} and show how the correspondence works in Section \ref{standardc} for standard contact structures and Section \ref{nonstandard} for the non-standard ones. We prove the theorem in Section \ref{proof} and show the transformation to $(2,3,5)$-distributions in Section \ref{dft}.

\section{Background material on contact structures of type $G_2$} \label{background}
We review the background theory of contact structures of type $G_2$ as given in \cite{ky93}, \cite{ky06} and \cite{LNS}. Let $J$ be a $5$-dimensional manifold and $\mathcal{C}=\{\varpi=0\}$ a linear differential system on $J$ on of codimension $1$. 
A Lie contact structure on $J$ is a 1-form $\varpi$ such that $\varpi \wedge \der \varpi \wedge \der\varpi$ is a volume form of $J$. 
By Darboux's theorem, there exists local coordinates $(x,y,z,p,q)$ on $J$ such that the contact 1-form $\varpi$ is given by
\[
\varpi=\der z-p \der x-q \der y.
\]
Let $(J, \mathcal{C})$ be such a Lie contact structure in dimension $5$.
Following \cite{ky93}, we can associate to the simple exceptional Lie group $G_2$ of non-compact type a Lie contact structure in the following manner. Let $G_2$ be the non-compact simple exceptional Lie group of dimension $14$ over $K=\R$ or $\C$ and $\frak{g}$ be its Lie algebra. We can associate to the Dynkin diagram $(G_2,\{\alpha_2\})$ on page 454 of \cite{ky93} the contact gradation
\[
\frak{g}=\frak{g}_{-2}\oplus\frak{g}_{-1}\oplus\frak{g}_{0}\oplus\frak{g}_{1}\oplus\frak{g}_{2}
\]
such that 
\[
[\frak{g}_p, \frak{g}_q]\in \frak{g}_{p+q} \qquad \mbox{for~} p, q \in \Z, 
\]
where $\dim{\frak{g}_{-2}}=1$ and $\frak{g}_{-1}\neq \{0\}$. Furthermore, we require that the Levi-bracket 
\[
[-,-]:\frak{g}_{-1}\times \frak{g}_{-1}\rightarrow \frak{g}_{-2} 
\]
is non-degenerate. 
Let 
\[
\frak{p}_2=\frak{g}_{0}\oplus\frak{g}_{1}\oplus\frak{g}_{2}
\]
be the parabolic subalgebra associated to this gradation.
There is also a gradation of depth $3$ associated to the Dynkin diagram 
$(G_2,\{\alpha_1\})$ with
\[
\frak{g}=\frak{g}_{-3}\oplus\frak{g}_{-2}\oplus\frak{g}_{-1}\oplus\frak{g}_{0}\oplus\frak{g}_{1}\oplus\frak{g}_{2}\oplus\frak{g}_{3}.
\]
Let 
\[
\frak{p}_1=\frak{g}_{0}\oplus\frak{g}_{1}\oplus\frak{g}_{2}\oplus \frak{g}_3
\]
be the parabolic subalgebra associated to this gradation.
Let $G=\rm{Int}(\frak{g})$ be the adjoint group of $\frak{g}$ and $R$ be the adjoint orbit of $G$ passing through an element $X_o$ in the graded piece $\frak{g}_{2}$. 
Let $P_1\cap P_2$ be the isotropy subgroup of $G$ that is the corresponding Borel subgroup of $G$ with Lie algebra $\frak{p}_1\cap \frak{p}_2$. Using the Killing form on $\frak{g}$, the covector $\om_o$ of $X_o$ is a left-invariant 1-form on $G$ equivariant with the principal right action of $P_1\cap P_2$ given by
\[
R_g^{*}\om_o={\rm Ad}^{*}(g^{-1}) \om_o=\om_o,
\]
where the group element $g \in P_1\cap P_2$ and $\om_o$ also annihilates the isotropic vectors. Thus $\om_o$ is projectable through $\pi:G\rightarrow R$ to a $G$-invariant 1-form, called $\al$, on $R=G/(P_1\cap P_2)$ with $\pi^{*}\al=\om_o$. The adjoint orbit $R$ then inherits a symplectic structure over $K^{\times}$ from the symplectic form $\der\alpha$ on $R$. 

Furthermore, let $\P(\frak{g})$ be the $K^{\times}$-projectivization of the vector fields of $\frak{g}$. Let $J$ be the $G$-orbit passing through $[X_o]= \frak{g}_2$ in $\P(\frak{g})$. There is a further projection $p:R\rightarrow J$ given by restricting the map
\[
p:\frak{g}\backslash \{0\} \rightarrow \P(\frak{g}).
\]
This gives $(R,J,p)$ the structure of a principal $K^{\times}$ bundle over $J$ with 
\[
R_g^{*}\alpha=k \cdot \alpha
\]
for $k\in K^{\times}$ and $g
\in P_2$. Furthermore, we have the splitting of $T(R)$ induced from 
$\al$ giving the vertical subbundle $\ker{p_{*}}\subset \ker \al=\{X\in T(R)|\alpha(X)=0\} \subset T(R)$, where $p_*$ is the pushforward map
\[
p_{*}:T(R)\rightarrow T(J).
\]
Let $\mathcal{C}(u)=p_{*}(\ker{\alpha}(x))$ be the Lagrangian subspace for the symplectic form $\der\al$ for each $u=p(x)\in J$. Then the $G$-invariant 1-form $\al$ on $R$ defines a $G$-invariant differential system $\mathcal{C}$ on $J$. In this sense, $(J, \mathcal{C})$ is called a Lie contact structure of type $G_2$. 

Starting from $(J, \mathcal{C})$, a Lie contact structure of type $G_2$ can be equivalently defined as a reduction of the structure group of the frame bundle of $\mathcal{C}$ to $GL(2,K)$ (2.7 of \cite{LNS}). In other words there is an isomorphism between $\mathcal{C}\cong S^3(V)$, where $V\rightarrow J$ is a rank $2$ bundle over $J$ and the Levi bracket is invariant under the action of $\frak{gl}(V)$ on $\mathcal{C}$. Since $GL(2,K)$ is the common stabilizer group of $\varpi$ and a symmetric $4$-form $\Up$, known as a structural tensor (4.46 of \cite{LNS}), a Lie contact structure of type $G_2$ on $J$ can therefore, for the purposes here, be completely described by the collection of 1-form and structural tensor $[\varpi, \Up]$.
 
For computation purposes in this paper later on, where we utilize the double fibration, it is better to consider the Lagrange-Grassman bundle 
\[
L(J)=\bigcup_{u\in J}L_u
\]
given by the union of the Lagrangian subspaces $L_u$ of $\mathcal{C}_u$ with the symplectic form $\der \varpi$.
There is a projection $\pi:L(J) \rightarrow J$. Let $U\subset L(J)$ be an open set in $L(J)$ that projects to $U^1\subset J$. Let $v\in U$ such that $\pi(v)=u \in U^1$. Then there exists coordinates such that
\begin{align*}
\der p|_{v}=-p_{11} \der x|_{v}-p_{12} \der y|_{v},\\
\der q|_{v}=-p_{21} \der x|_{v}-p_{22} \der y|_{v},
\end{align*} 
where the compatibility condition $p_{12}=p_{21}$ is satisfied. (Our choice of $p_{11}$, $p_{12}$, $p_{22}$ is different from those in \cite{ky06} by a sign).  
Denoting by 
\begin{align*}
\om_2&=\der p+p_{11}\der x+p_{12}\der y,\\
\om_3&=\der q+p_{21}\der x+p_{22}\der y,\\
\om_4&=\hspace{24pt} p_{31} \der x+p_{32}\der y,
\end{align*}
the projection $\pi$ induces the map
\begin{center}
\begin{tikzcd}[column sep=small]
T(L(J))
  \arrow[d,"\pi_{*}"]\arrow[r,phantom,"\supset" description] 
 &E(u)=\pi_{*}^{-1}(u)  \arrow[d,"\pi_{*}"]
\\
T(J)\arrow[r,phantom,"\supset" description] &u
\end{tikzcd}
\end{center}
where $\pi_{*}^{-1}(u)=E(u)=\{\varpi=\om_2=\om_3=\om_4=0\}$. The functions $(x,y,z,p,q,p_{11},p_{12},p_{22},p_{31},p_{32})$ define local coordinates on $L(J)$ and we pass to the submanifold $R \subset L(J)$ such that $R\rightarrow J$ is a submersion. 
We now take further
\begin{align*}
R=\big\{p_{11}=t^2H_t-2tH+2\int H\der t,\hspace{12pt}p_{12}=H-tH_t,\hspace{12pt}p_{22}=H_t,\hspace{12pt}p_{31}=-t,\hspace{12pt}p_{32}=1\big\},
\end{align*}
where $H(t)$ is a single function of the fibre variable $t$. 
The general form of these Lie contact 1-forms that give rise to $G_2$ symmetry are of the form 
\begin{align*}
\om_2&=\der p+(t^2H_t-2tH+2\int H\der t)\der x+(H-tH_t)\der y,\\
\om_3&=\der q+(H-tH_t)\der x+H_t\der y,\\
\om_4&=\der y-t \der x,
\end{align*}
where $H(t)$ satisfies Noth's equation
\begin{equation}\label{h-noth}
10(H'')^3H^{(6)}-70(H'')^2H^{(3)}H^{(5)}-49(H'')^2(H^{(4)})^2+
280H''(H^{(3)})^2H^{(4)}-175(H^{(3)})^4=0.
\end{equation}
Here the derivative taken with respect to $t$ variable. This is the computation on page 180 of \cite{Cartan1910}.
We first review the case where the solution of Noth's equation is taken to be $H(t)=3t^2$, which gives rise to the standard contact structure of type $G_2$ associated to the twisted cubic and Proposition 4.4 in \cite{ufo1}.

\section{Contact structures of type $G_2$ associated to the solution $H(t)=3t^2$ of Noth's equation}\label{standardc}

When $H(t)=3t^2$, Noth's equation is satisfied and the structural 1-forms up to translation in the $p_{11}$ term are 
\begin{align*}
\om_2&=\der p+2t^3\der x-3t^2\der y,\\
\om_3&=\der q-3t^2\der x+6t\der y,\\
\om_4&=\der y-t \der x.
\end{align*}
We now consider $\om_2$, $\om_3$ and $\om_4$ as polynomials in $t$ with coefficients in the ring of differentials $\R[dx,dy,dz,dp,dq]$ where multiplication as symbols is formal. Let us call
\begin{align*}
p_2&=2 dx t^3-3 dy t^2+dp,\\
p_3&=-3 dx t^2+6 dy t+dq,\\
p_4&=-dx t+dy.
\end{align*}

The condition that these three polynomials share a common zero means that the resultant of any two of them is zero. We find
\begin{align*}
{\rm Res}(p_2, p_3)&=-dx(27dp^2dx^2+54dpdqdxdy+108dpdy^3-4dq^3dx-9dq^2dy^2),\\ 
{\rm Res}(p_3, p_4)&=-dx(dx^2dp-dy^3),\\ 
{\rm Res}(p_2, p_4)&=-dx(dq dx+3 dy^2).
\end{align*}
The vanishing of these coefficients descend to give the vanishing of the symmetric $2$, $3$ and $4$-tensors
\begin{align*}
\Up&=27\der p^2\der x^2+54\der p\der q \der x \der y+108\der p\der y^3-4\der q^3\der x-9\der q^2\der y^2&\in \Ga({\rm Sym}^4(T^*J)),\\ 
\mu&=\der x^2\der p-\der y^3&\in \Ga({\rm Sym}^3(T^*J)),\\ 
g_1&=\der q \der x+3 \der y^2&\in \Ga({\rm Sym}^2(T^*J)).
\end{align*}
In the literature \cite{ufo1}, $\Up$ is called a structural tensor and together with $\varpi$ is used to determine a Lie contact structure of type $G_2$ completely. 
Let
\begin{align*}
g_2=9 \der p \der y-\der q^2,\quad g_3=3 \der p \der x+\der y \der q. 
\end{align*}
Then 
\[
\Up=4 g_1 g_2+3 g_3^2 \quad \mbox{~and~}\quad  \mu=\frac{1}{3}(\der x g_3-\der y g_1). 
\]
The tensors $[\varpi, \Up]$ determine a Lie contact structure of $G_2$ type associated to the solution $H(t)=3t^2$ of Noth's equation and $(g_1,g_2,g_3)$ form a $GL(2,K)$-invariant module of symmetric bilinear forms. The full symmetry algebra is given in Proposition 4.4 of \cite{ufo1}. We shall call such a contact structure associated to the solution $H(t)=3t^2$ standard. In \cite{ufo1}, it is shown that given a Lie contact structure of type $G_2$ with reduction of the structure group of the frame bundle $\mathcal{C}$ to $GL(2,K)$, one can find coordinates associated to the twisted cubic to express the Lie contact structure as the standard one (Proposition 4.2 of \cite{ufo1}). Since standard contact structures of type $G_2$ are related to maximally symmetric $(2,3,5)$-distribution associated to the $1$-forms
\begin{align}\label{hc1}
\der Y-P \der X, \quad \der P-Q \der X \quad\mbox{and} \quad \der Z-Q^2\der X
\end{align}
by the double fibration transform, this result of choosing coordinates to express the Lie contact structure as a standard one can be seen as the dual version of the result of choosing coordinates such that the maximally symmetric $(2,3,5)$-distribution is given by the Hilbert-Cartan distribution (\ref{hc1}). With the exception of \cite{Cartan1910}, all the literature on contact strucutres of type $G_2$ up to now just treats the standard case. However, it is already known in \cite{Cartan1910} that the general moduli of contact structures of type $G_2$ is parametrized by solutions to Noth's equation, and the attention in the paper is focused on these other solutions obtained in \cite{r17}. Let us call the contact structures associated to a solution of Noth's equation that is not $SL(2,K)$ equivalent to $H(t)=3t^2$ non-standard. In this paper we are interested in the non-standard contact structure associated to the parametrization of the solution of Noth's equation obtained in Theorem 6.7 of \cite{r17}. We will also show in Section \ref{dft} later on, that the non-standard contact structures associated to the solutions of Noth's equation are dual to $(2,3,5)$-distributions of the form
\begin{align*}
\der Y-P \der X, \quad \der P-Q \der X, \quad \der Z-\frac{Q^2}{H_{XX}}\der X,
\end{align*}
as studied in \cite{r19} and \cite{r21c}, where again $H(X)$ satisfies Noth's equation.

 \section{Contact structures of type $G_2$ associated to solutions of Noth's equation}\label{nonstandard}
Noth's equation can be reduced to the generalized Chazy equation with parameter $\frac{3}{2}$ and its solution as parametrized by Schwarz triangle functions is presented in \cite{r17}. A rational parametrization of the solution is obtained in Theorem 6.7 of \cite{r17} and we consider this solution here. The equation is also dual to an equation that can be reduced to the generalized Chazy equation with parameter $\frac{2}{3}$. We will not discuss this situation in the paper. 
We now consider the case where the solution of Noth's equation is taken to be that given in Theorem 6.7 of \cite{r17}. The solutions are given by the curve
\begin{align}\label{noth1}
(t,H)=\left(\frac{2s^{\frac{1}{3}}}{s+2},-\frac{4s^{\frac{2}{3}}}{s+2}\right)
\end{align}
and
\begin{align}\label{noth2}
(t,H)=\left(\frac{2s^{\frac{2}{3}}}{1+2s},-\frac{4s^{\frac{1}{3}}}{1+2s}\right)
\end{align}
obtained from inverting $s \mapsto \frac{1}{s}$ in (\ref{noth1}).
The first parametrization (\ref{noth1}) gives the structural 1-forms
\begin{align*}
\om_2&=\der p+\frac{4}{3(s-1)}\der x-\frac{2s^{\frac{2}{3}}}{s-1}\der y,\\
\om_3&=\der q-\frac{2s^{\frac{2}{3}}}{s-1}\der x-\frac{s^{\frac{1}{3}}(s-4)}{s-1}\der y,\\
\om_4&=\der y-\frac{2s^{\frac{1}{3}}}{s+2}\der x.
\end{align*}
Substituting $r=s^{\frac{1}{3}}$ and eliminating the parameter $r$ by taking resultants implies the vanishing of the tensors
\begin{align}\label{upa}
\Up&=81\der p^4 - 216\der p^3\der x + 216\der p^2\der q\der y + 144\der p^2\der x^2+24\der p\der q^3 \\  \nonumber
&\quad- 288\der p\der q\der x\der y - 384\der p\der y^3 - 48\der q^2\der y^2 + 512\der x\der y^3,\\\label{ma1}
\mu_1&=27\der p^3-36\der p^2\der x+32\der y^3,\\\label{ma2}
\mu_2&=\der q^3-8\der q\der x\der y + 16\der y^3,\\
\nu&=4\der x^2 + 6\der x\der y + 9\der y^2 \mbox{~and~}\ka=2\der x-3 \der y\nonumber
\end{align}
on the locus $\{\om_2=\om_3=\om_4=0\}$.
Let us consider the module of symmetric bilinear forms given by the span of
\begin{align*}
g_1&=8\der y^2-3\der p \der q,\\
g_2&=9 \der p^2-12\der x \der p+4\der y \der q,\\
g_3&=\der q^2+6 \der y \der p-8 \der x \der y.
\end{align*}
The structural tensors $\Up$ given in (\ref{upa}) and $\mu_1$, $\mu_2$ given in (\ref{ma1}), (\ref{ma2}) have the property that
\begin{align*}
\mathcal{L}_{X}\Up&=0 \mod \varpi, \Up,\\
\mathcal{L}_{X}\mu_1&=0 \mod \varpi, g_1, g_2, g_3,\\
\mathcal{L}_{X}\mu_2&=0 \mod \varpi, g_1, g_2, g_3.
\end{align*}

Here the $X$'s are the vector fields in the symmetry algebra of the 1-forms given by the exterior differential system $\{\varpi=\om_2=\om_3=\om_4=0\}$. 
We also have
\begin{align*}
\Up=g_2^2-8 g_1g_3, \quad \mu_1=4 \der y g_1+\der p g_2\mbox{~and~}\mu_2=2 \der y g_1+\der q g_3.
\end{align*}
It is unclear what role does the 1-form $\ka$ and the 2-form $\nu$ play here. 

\begin{theorem}\label{thm1}
The symmetry algebra of the contact structure of $G_2$ type with the contact form $\varpi$ and structural tensor $\Up$ given by (\ref{upa}), associated to the solution of Noth's equation with the solution parametrized by (\ref{noth1}) is split $\frak{g}_2$ with the vector fields given by
\begin{align*}
S^1&=-(-3qx-6y^2+\frac{9}{4}pq)\partial_x-(3z-3px-4qy+\frac{9}{8}p^2)\partial_y-(-8xy^2-3pqx-2q^2y+\frac{9}{4}p^2q)\partial_z\\&\quad{}+8y^2\partial_p-(12py-q^2-16xy)\partial_q,\\
L^2&=-(-\frac{9}{2}px+\frac{27}{16}p^2+\frac{3}{2}z+2x^2)\partial_x-\frac{1}{4}q^2\partial_y-(-\frac{9}{4}p^2x+\frac{9}{8}p^3+\frac{1}{6}q^3+2x z)\partial_z\\&\quad{}-(2z-2px+\frac{3}{4}p^2)\partial_p+(\frac{3}{2} pq-2 qx)\partial_q,\\
S^3&=-\frac{3}{2}y\partial_x-\frac{1}{2}q\partial_y-(2 x y+\frac{1}{4}q^2)\partial_z-2y \partial_p+(\frac{3}{2}p-2x)\partial_q,\\
L^{4}&=\frac{3}{4}\partial_x+x\partial_z+\partial_p,\quad S^{5}=\partial_y,\quad L^{6}=\partial_z,\quad S^{7}=y\partial_z+\partial_q,\quad L^{8}=\partial_{x},\\
S^{9}&=\frac{3}{4}q\partial_x+\frac{3}{4}p\partial_y+(\frac{3}{4}pq+2y^2)\partial_z+4y \partial_q,\\
L^{10}&=(\frac{3}{4}qy+\frac{27}{32}p^2-\frac{3}{4}z)\partial_x+\frac{3}{4}py\partial_y+(\frac{3}{4}pqy+\frac{9}{16}p^3+\frac{2}{3}y^3)\partial_z+\frac{3}{4}p^2\partial_p+2y^2\partial_q,\\
S^{11}&=(\frac{3}{8}q^2+\frac{9}{4}py)\partial_x+(\frac{3}{4}pq+y^2)\partial_y+(\frac{3}{4}pq^2+\frac{9}{8}p^2y+3yz)\partial_z+3py\partial_p-(-3z+\frac{9}{8}p^2-q y)\partial_q,\\
L^{12}&=(3zx+\frac{1}{8}q^3-3qxy - 2y^3 - \frac{27}{8}p^2x+\frac{27}{16}p^3+\frac{9}{4}pqy)\partial_x+(3yz - 3pxy +\frac{3}{8}pq^2- 2qy^2 + \frac{9}{8}p^2y)\partial_y\\
&\quad{}+(-\frac{9}{4}p^3x-\frac{8}{3}xy^3+3z^2+\frac{3}{8}q^3p-q^2y^2+\frac{9}{4}qy p^2+\frac{81}{64}p^4-3 x p q y)\partial_z\\&\quad{}+(3pz-3p^2 x-\frac{8}{3}y^3+\frac{9}{8}p^3)\partial_p-(-3 q z+\frac{9}{8}p^2 q-6 p y^2+q^2 y+8 x y^2)\partial_q,\\
h_1&=-(\frac{9}{4}p-\frac{3}{2}x)\partial_x-\frac{1}{2}y\partial_y-\frac{9}{8}p^2\partial_z-\frac{3}{2}p \partial_p+\frac{1}{2}q\partial_q,\\
h_2&=\frac{\sqrt{3}}{2}(x\partial_x+y\partial_y+2 z \partial_z+p \partial_p+q\partial_q).
\end{align*}
\end{theorem}
 With this choice of the Cartan subalgebra spanned by $h_1$ and $h_2$, the root diagram is given by the picture below. 
\begin{figure}[h!]
\begin{tikzpicture}
	\draw [stealth-stealth](-1,0) -- (1,0);
\draw (1,0) node[anchor=west] {{\tiny $S^3$}};
\draw (-1,0) node[anchor=east] {{\tiny $S^9$}};
	\draw [stealth-stealth](0,-1.732) -- (0,1.732);
\draw (0,-1.732) node[anchor=north] {{\tiny $L^6$}};
\draw (0,1.732) node[anchor=south] {{\tiny $L^{12}$}};
\draw [stealth-stealth](-0.5,-0.866) -- (0.5,0.866);
\draw (-0.5,-0.866) node[anchor=north] {{\tiny $S^{7}$}};
\draw (0.5,0.866) node[anchor=south] {{\tiny $S^1$}};
\draw [stealth-stealth](-1.5,-0.866) -- (1.5,0.866);
\draw (-1.5,-0.866) node[anchor=north] {{\tiny $L^8$}};
\draw (1.5,0.866) node[anchor=south] {{\tiny $L^2$}};
\draw [stealth-stealth](1.5,-0.866) -- (-1.5,0.866);
\draw(1.5,-0.866) node[anchor=north] {{\tiny $L^4$}};
\draw (-1.5,0.866) node[anchor=south] {{\tiny $L^{10}$}};
\draw [stealth-stealth](0.5,-0.866) -- (-0.5,0.866);
\draw (0.5,-0.866) node[anchor=north] {{\tiny $S^{5}$}};
\draw (-0.5,0.866) node[anchor=south] {{\tiny $S^{11}$}};
\end{tikzpicture}
\end{figure} 
We have labelled the vectors in the diagram accordingly so that each vector is either $S^i$ or $L^i$, where $S$ and $L$ depends on whether the root is short or long and $i$ denotes the hour on a standard clock when superimposing the root diagram on a clock face.

The coordinate diffeomorphism to the standard contact structure of type $G_2$ is given as follows. Consider the map
\begin{align*}
(x_1,y_1,p_1,q_1,z_1)=\left(\frac{1}{4}(3x-p),-\frac{12^{2/3}}{8}y,x,-\frac{12^{1/3}}{6}q,\frac{1}{4}(z-xp+\frac{3}{2}x^2)\right).
\end{align*}
Then we find
\begin{align*}
\der z_1-p_1\der x_1-q_1\der y_1=\frac{1}{4}\varpi,\\
8\der y_1\der y_1-3\der p_1\der q_1=\frac{12^{1/3}}{2}g_1,\\
9\der p_1\der p_1-12 \der p_1\der x_1+4 \der y_1\der q_1=g_3,\\
\der q_1\der q_1+6 \der y_1 \der p_1-8 \der x_1\der y_1=-\frac{12^{2/3}}{36}g_2,
\end{align*}
where $g_1$, $g_2$ and $g_3$ are standard symmetric bilinear forms given in Section \ref{standardc}. This coordinate diffeomorphism induces the isomorphism of $S^3(V)$ on the standard contact structure with $\mathcal{C}$. 

Conversely, given such a coordinate diffeomorphism, we can find a solution to Noth's equation as follows. To get back the original parametrization, we need to apply a translation in the $p_{11}$ term.
Let $\varphi:J\rightarrow J_{standard}$ be the coordinate diffeomorphism. Let $t=r$ be the fibre parameter of the projection $\pi_{2}:R\rightarrow J$. Then the 1-forms $\varpi$ together with $\om_2$, $\om_3$, $\om_4$ given by
\begin{align*}
\om_2&=\der p+(2r^3-5)\der x-3r^2\der y,\\
\om_3&=\der q-3r^2\der x+6r\der y,\\
\om_4&=\der y-r \der x,
\end{align*}
on $J_{standard}$ pulls back to the 1-forms
\begin{align*}
\der p_1+p_{11}(r) \der x_1+p_{12}(r)\der y_1,\\
\der q_1+p_{21}(r) \der x_1+p_{22}(r)\der y_1,\\
\der z_1-p_1\der x_1-q_1\der y_1,\\
\der y_1+p_{31}(r)\der x_1
\end{align*}
on $J$ for some functions $p_{11}(r)$, $p_{12}(r)=p_{21}(r)$, $p_{22}(r)$ and $p_{31}(r)$.
These 1-forms must lie in the span of $\om_2$, $\om_3$, $\om_4$, from which we can solve for the functions to obtain
\begin{align*}
p_{11}(r)&=-\frac{2}{r^3-1},\quad p_{12}(r)=12^{1/3}\frac{r^2}{r^3-1},\\
p_{22}(r)&=12^{2/3}\frac{r(r^3-4)}{6(r^3-1)},\quad p_{31}(r)=-\frac{12^{2/3}r}{2(r^3+2)},
\end{align*}
so that
\begin{align*}
\der z_1-p_1\der x_1-q_1\der y_1&=\frac{1}{4}\varpi,\\
\der p_1-\frac{2}{r^3-1} \der x_1+12^{1/3}\frac{r^2}{r^3-1}\der y_1&=\frac{1}{2(r^3-1)}\om_2,\\
\der q_1+12^{1/3}\frac{r^2}{r^3-1} \der x_1+12^{2/3}\frac{r(r^3-4)}{6(r^3-1)}\der y_1&=-\frac{12^{1/3}r^2}{4(r^3-1)}\om_2-\frac{12^{1/3}}{6}\om_3,\\
\der y_1-\frac{12^{2/3}r}{2(r^3+2)}\der x_1&=\frac{12^{2/3}r}{8(r^3+2)}\om_2+\frac{12^{2/3}(r^3-1)}{4(r^3+2)}\om_4.
\end{align*}
This gives the parametrization 
\begin{align*}
(t,H)=\left(\frac{12^{2/3}r}{2(r^3+2)},2\cdot 12^{1/3}\frac{r^2}{r^3+2}\right)
\end{align*}
of the solution to (\ref{h-noth}).
If we used the original standard $1$-form $\om_2$, not translating the $p_{11}$ term, we would obtain a different parametrization
\begin{align*}
(t,H)=\left(\frac{12^{2/3}r}{2(r^3-3)},2\cdot 12^{1/3}\frac{r^2}{r^3-3}\right)
\end{align*}
of the solution to (\ref{h-noth}).

We now consider the parametrization of Noth's equation given by 
\begin{align*}
(t,H)=\left(\frac{2s^{\frac{2}{3}}}{1+2s},-\frac{4s^{\frac{1}{3}}}{1+2s}\right)
\end{align*}
obtained from inverting $s \mapsto \frac{1}{s}$. Although inversion of $s$ is a fractional linear transformation, this is not an $SL(2,K)$-equivalent solution to the one given by (\ref{noth1}) in the sense of \cite{r17}, as we are changing the Schwarz triangle function parametrizing the solution. In the case of inversion, this corresponds to a permutation of the angles in the corresponding Schwarz triangle function. 
This solution gives the structural 1-forms
\begin{align*}
\om_2&=\der p-\frac{4}{3(s-1)}\der x+\frac{2s^{\frac{1}{3}}}{s-1}\der y,\\
\om_3&=\der q+\frac{2s^{\frac{1}{3}}}{s-1}\der x+\frac{1-4s}{s^{\frac{1}{3}}(s-1)}\der y,\\
\om_4&=\der y-\frac{2s^{\frac{2}{3}}}{2s+1}\der x.
\end{align*}
Under the inversion $s \mapsto \frac{1}{s}$, $\om_3$ and $\om_4$ agrees with that given above but we have
\[
\om_2=\der p-\frac{4}{3(1/s-1)}\der x+\frac{2}{s^{\frac{1}{3}}(1/s-1)}\der y=\der p+\frac{4s}{3(s-1)}\der x-\frac{2s^{\frac{2}{3}}}{s-1}\der y
\]
which does not agree with $\om_2$ above. 
This is because the fractional linear transformation $s\mapsto \frac{1}{s}$ does not leave the $p_{11}$ term given by
\[
t^2H_t-2tH+2\int H\der t
\] invariant. 
 Eliminating the parameter $s$ gives the tensors
\begin{align}\label{upb}
\Up&=81\der p^4+216\der p^3\der x + 216\der p^2\der q\der y + 144\der p^2\der x^2+24\der p\der q^3 \\ \nonumber
&\quad+288\der p\der q\der x\der y - 384\der p\der y^3 +32\der q^3\der x-48\der q^2\der y^2,\\\label{mb1}
\mu_1&=27 \der p^3+72 \der p^2 \der x+48 \der p \der x^2+32 \der y^3,\\\label{mb2}
\mu_2&=\der q^3-8\der q\der x\der y + 16\der y^3,\\ \nonumber
\nu&=4\der x^2 + 6\der x\der y + 9\der y^2,\\
\ka&=2\der x-3\der y.\nonumber
\end{align}

Let us form the module of symmetric bilinear forms by taking
\begin{align*}
g_1&=8\der y^2-3\der p \der q-4\der x \der q,\\
g_2&=9 \der p^2+12\der x \der p+4\der y \der q,\\
g_3&=\der q^2+6 \der y \der p.
\end{align*}
We have
\begin{align*}
\Up=g_2^2-8 g_1g_3, \quad \mu_1=4\der y g_1+(4\der x+3 \der p)g_2\mbox{~and~}\mu_2=2\der y g_1+\der q g_3.
\end{align*}
Furthermore, 
\begin{align*}
\mathcal{L}_{X}\Up&=0 \mod \varpi, \Up,\\
\mathcal{L}_{X}\mu_1&=0 \mod \varpi, g_1, g_2, g_3,\\
\mathcal{L}_{X}\mu_2&=0 \mod \varpi, g_1, g_2, g_3,
\end{align*}
where $X$ is a vector field in the symmetry algebra of the differential system $\{\varpi=\Up=0\}$. 
\begin{theorem}\label{thm2}
The symmetry algebra of the contact structure of $G_2$ type with the contact form $\varpi$ and structural tensor $\Up$ given by (\ref{upb}) associated to the solution of Noth's equation with the parametrization given by (\ref{noth2}) is split $\frak{g}_2$ with the vector fields given by
\begin{align*}
S^1&=-(\frac{3}{16}q^2+\frac{3}{2}xy+\frac{9}{8}py)\partial_x-(\frac{3}{8}pq+\frac{1}{2}y^2+\frac{1}{2}qx)\partial_y-(\frac{9}{16}p^2y+\frac{1}{4}q^2x+\frac{3}{8}q^2p+\frac{3}{2}yz)\partial_z\\&\quad{}+\frac{1}{4}q^2\partial_p+(-\frac{3}{2}z+\frac{3}{2}px-\frac{1}{2}qy+\frac{9}{16}p^2)\partial_q,\\
L^2&=-(\frac{9}{4}px+\frac{27}{32}p^2-\frac{3}{4}z+x^2+\frac{3}{4}qy)\partial_x-(\frac{3}{4}py+xy)\partial_y-(\frac{9}{8}p^2x+\frac{9}{16}p^3+x z+\frac{3}{4}pqy+\frac{2}{3}y^3)\partial_z\\&\quad{}-(z-px-qy-\frac{3}{8}p^2)\partial_p-2y^2\partial_q,\\
S^3&=-\frac{3}{8}q\partial_x-(\frac{3}{8}p+\frac{1}{2}x)\partial_y-(\frac{3}{8}pq+y^2)\partial_z+\frac{1}{2}q\partial_p-2y \partial_q,\\
L^{4}&=-\frac{3}{4}\partial_x+x\partial_z+\partial_p,\quad S^{5}=y\partial_z+\partial_q,\quad L^{6}=\partial_z,\quad S^{7}=\partial_y,\quad L^{8}=\partial_{x},\\
S^9&=3y\partial_x+q\partial_y+\frac{1}{2}q^2\partial_z-3p\partial_q,\\
L^{10}&=(3z+\frac{27}{8}p^2)\partial_x+\frac{1}{2}q^2\partial_y+(\frac{9}{4}p^3+\frac{1}{3}q^3)\partial_z-3p^2\partial_p-3pq\partial_q,\\
S^{11}&=(-6y^2+\frac{9}{4}qp)\partial_x+(3z+\frac{9}{8}p^2-4qy)\partial_y+(\frac{9}{4}p^2q-2yq^2)\partial_z-3pq\partial_p+(12py-q^2)\partial_q,\\
L^{12}&=-(3zx+\frac{1}{8}q^3+\frac{9}{4}pqy - 2y^3+\frac{27}{8}p^2x+\frac{27}{16}p^3)\partial_x-(3yz +\frac{1}{2}q^2x +\frac{3}{8}pq^2- 2qy^2 + \frac{9}{8}p^2y)\partial_y\\
&\quad{}-(\frac{9}{4}p^3x-\frac{1}{3}xq^3+3z^2+\frac{3}{8}q^3p-q^2y^2+\frac{9}{4}qy p^2+\frac{81}{64}p^4-3 x p q y)\partial_z\\&\quad{}+(-3pz+3p^2 x+3pqy+\frac{1}{6}q^3+\frac{9}{8}p^3)\partial_p+(-3 q z+3pqx+\frac{9}{8}p^2 q-6 p y^2+q^2 y)\partial_q,\\
h^1&=(\frac{9}{4}p+\frac{3}{2}x)\partial_x+\frac{1}{2}y\partial_y+\frac{9}{8}p^2\partial_z-\frac{3}{2}p \partial_p-\frac{1}{2}q\partial_q,\\
h^2&=\frac{\sqrt{3}}{2}(x\partial_x+y\partial_y+2 z \partial_z+p \partial_p+q\partial_q).
\end{align*}
Again with this choice of the Cartan subalgebra spanned by $h_1$ and $h_2$, the root diagram is given by the picture above in clock face notation. 
\end{theorem}

The coordinate diffeomorphism to the standard contact structure of type $G_2$ is given as follows. Consider
\begin{align*}
(x_2,y_2,p_2,q_2,z_2)=\left(\frac{1}{4}(3x+p),\frac{12^{2/3}}{8}y,-\frac{1}{3}p,-\frac{12^{1/3}}{6}q,-\frac{1}{4}(z+\frac{1}{6}p^2)\right).
\end{align*}
Then we find
\begin{align*}
\der z_2-p_2\der x_2-q_2\der y_2=-\frac{1}{4}\varpi,\\
8\der y_2\der y_2-3\der p_2\der q_2-4\der x_2\der q_2=\frac{12^{1/3}}{2}g_1,\\
9\der p_2\der p_2+12 \der p_2\der x_2+4 \der y_2\der q_2=-g_3,\\
\der q_2\der q_2+6 \der y_2 \der p_2=-\frac{12^{2/3}}{36}g_2,
\end{align*}
where $g_1$, $g_2$ and $g_3$ are again the standard symmetric bilinear forms in Section \ref{standardc}. This coordinate diffeomorphism gives the isomorphism of $\mathcal{C}$ with $S^3(V)$ of the standard contact structure. 
There is a natural $GL(2,K)$ action on $J_{standard}$ which rescales the 1-forms and hence the structural tensor associated to $\mathcal{C}_{standard}$ by a constant. 
Conversely, given such a coordinate diffeomorphism, we can find a solution to Noth's equation as follows. 
Let $\varphi:J\rightarrow J_{standard}$ be the coordinate diffeomorphism. It turns out that in order to obtain the solution given by the parametrization (\ref{noth2}), we should apply the same transformation $r \to \frac{1}{r}$ on the 1-forms on $(J_{standard}, \mathcal{C}_{standard})$ followed by an additonal translation of the $p_{11}$ term in $\om_2$ to get the 1-forms 
$\om_2$, $\om_3$, $\om_4$ given by
\begin{align*}
\om_2&=\der p+(\frac{2}{r^3}+1)\der x-\frac{3}{r^2}\der y,\\
\om_3&=\der q-\frac{3}{r^2}\der x+\frac{6}{r}\der y,\\
\om_4&=\der y-\frac{1}{r}\der x.
\end{align*}
These 1-forms on $J_{standard}$ pulls back to the 1-forms
\begin{align*}
\der z_2-p_2\der x_2-q_2\der y_2,\qquad& \der p_2+p_{11}(r) \der x_2+p_{12}(r)\der y_2,\\
\der y_2+p_{31}(r)\der x_2,\qquad& \der q_2+p_{12}(r) \der x_2+p_{22}(r)\der y_2,
\end{align*}
on $J$ for some functions $p_{11}(r)$, $p_{12}(r)$, $p_{22}(r)$ and $p_{31}(r)$.
These 1-forms must lie in the span of $\om_2$, $\om_3$, $\om_4$, which implies that 
\begin{align*}
p_{11}(r)&=-\frac{2(r^3+2)}{3(r^3-1)},\quad p_{12}(r)=p_{21}(r)=12^{1/3}\frac{r}{r^3-1},\\
p_{22}(r)&=-12^{2/3}\frac{4r^3-1}{6r(r^3-1)},\quad p_{31}(r)=-\frac{12^{2/3}r^2}{4r^3+2},
\end{align*}
so we have
\begin{align*}
\der p_2-\frac{2(r^3+2)}{3(r^3-1)} \der x_2+12^{1/3}\frac{r}{r^3-1}\der y_2&=-\frac{r^3}{2r^3-2}\om_2,\\
\der q_2+12^{1/3}\frac{r}{r^3-1} \der x_2-12^{2/3}\frac{4r^3-1}{6r(r^3-1)}\der y_2&=\frac{12^{1/3}r}{4(r^3-1)}\om_2-\frac{12^{1/3}}{6}\om_3,\\
\der z_2-p_2\der x_2-q_2\der y_2&=-\frac{1}{4}\varpi,\\
\der y_2-\frac{12^{2/3}r^2}{4r^3+2}\der x_2&=-\frac{12^{2/3}r^2}{8(2r^3+1)}\om_2+\frac{12^{2/3}(r^3-1)}{4(2r^3+1)}\om_4.
\end{align*}
This gives the parametrization 
\begin{align*}
\left(\frac{12^{2/3}r^2}{4r^3+2},-\frac{2\cdot 12^{1/3}r}{2r^3+1}\right)
\end{align*}
of the solution to Noth's equation. 
If we did not apply the translation in the $p_{11}$ term, we would have obtained the solution parametrized by 
\begin{align*}
\left(\frac{12^{2/3}r^2}{2(3r^3+1)},-\frac{2\cdot 12^{1/3}r}{3r^3+1}\right)
\end{align*}
instead. We therefore see that the correspondence requires us to apply a translation in the $p_{11}$ term and also to fix the fibre parameter of the projection given by $\pi_2$. 

\section{Proof of Theorem \ref{maintheorem}} \label{proof}

The proof is derived from observing the above correspondence. In principle, the theorem should be proven in the most general way instead of relying on certain parametrization of solutions of Noth's equation. Nonetheless from the above solutions in Section \ref{nonstandard}, we constructed the local coordinate diffeomorphism from $(J,\mathcal{C})$ to the standard contact structure of type $G_2$. Conversely, starting from such a local coordinate diffeomorphism, we also showed how to construct a solution to Noth's equation. 
\begin{proof}
Let $(J,\mathcal{C})$ be an abstract contact structure of type $G_2$ and $(J_{standard}, \mathcal{C}_{standard})$ be the standard one. 
Given any solution to Noth's equation, the exterior differential system (E.D.S.) associated to $\tilde D_2$ has vanishing Cartan invariant and so $J$ has the local structure of $G_2/P_2$ with maximal symmetry. By uniqueness of the twisted cubic standard contact structure, up to local diffeomorphism there exists $\varphi$ such that
\[
\varphi:\mathcal{C}\rightarrow \mathcal{C}_{standard}=S^3(V)
\]
is an isomorphism. 
Thus $\varphi$ gives the local coordinate diffeomorphism and we can pull-back the contact structure $\mathcal{C}_{standard}$ to $\mathcal{C}$ on $J$, where the fibre parameter is the natural independent parameter in the solution to Noth's eqaution. 
What is perhaps more interesting is the converse. 

Let $\varphi:\mathcal{C}\rightarrow S^3(V)$ be the given isomorphism that identifies an abstract $(J,\mathcal{C})$ with $(J_{standard}, \mathcal{C}_{standard})$ and $\mathcal{C}_{standard}=S^3(V)$. Then at the level of differential forms, it induces $\varphi^*$ such that the pulled back differential system lies in the span of the differential system of the standard one, that is to say 
\[
\varphi^*(\varpi, \om_2,\om_3,\om_4)\subseteq (\varpi, \om_2,\om_3, \om_4),
\]
where $\om_2$ is defined up to translation in the $p_{11}$ term. 
This means there exists coordinates $(x_1,y_1,p_1,q_1,z_1)$ such that
\begin{align*}
dz_1&=-p_1dx_1-q_1dy_1,\\
dp_1&=-p_{11}dx_1-p_{12}dy_1,\\
dq_1&=-p_{21}dx_1-p_{22}dy_1,\\
dy_1&=-p_{31}dx_1,
\end{align*}
where $p_{11}$, $p_{12}$, $p_{22}$ and $p_{31}$ are functions of the fibre parameter $r$ chosen in the projection map of $\pi_2$ onto $J_{standard}$. Once we fix this parameter $r$ in the projection to the standard contact structure, the independent variable in Noth's equation becomes a function of it. We now identify the differential system with the Lie 1-forms in Section \ref{background}. In this case, the pair $(-p_{31},-\int p_{22}\der p_{31})$ gives an identification with the solution $(t, H)$ of Noth's equation. 
\end{proof}
The identity morphism 
$\varphi=\rm{Id}$ corrsponds to the standard contact structure with solution of Noth's equation given by $3t^2$.

\section{Transformation to $(2,3,5)$-distributions}\label{dft}
We now discuss the double fibration transform and see how to get the corresponding $(2,3,5)$-distribution from the contact structure discussed above. 
It is known through Section 4.5 of \cite{ufo1} and \cite{ky06} that the double fibration
\begin{center}
\begin{tikzcd}
&  
G_2/(P_{1}\cap P_2)
  \arrow[dr,"\pi_2"]
  \arrow[dl,"\pi_1",swap] 
\\
G_2/P_1
&&
G_{2}/P_{2}
\end{tikzcd}
\end{center}
can be expressed using the notation of linear differential systems as:
\begin{center}
\begin{tikzcd}
&  
(R, \tilde D_1, \tilde D_2)
  \arrow[dr,"\pi_2"]
  \arrow[dl,"\pi_1",swap] 
\\
(M,D_1)
&&
(J,D_2)
\end{tikzcd}
\end{center}
The 5-dimensional manifold $M$ is equipped with the structure of a  $(2,3,5)$-distribution expressed through the linear differential system $D_1$, which can be encoded by the 1-forms
\begin{align*}
o_1=\der Y-P\der X,\quad o_2=\der P-Q\der X,\quad o_3=\der Z-\frac{Q^2}{H_{XX}}\der X,
\end{align*} with $D_1=\{o_1=o_2=o_3=0\}$.
The differential system on $M$ can be pulled back to that on $R$  by augmenting the 1-form 
\begin{align*}
o_4=\der Q-L\der X, 
\end{align*}
where $L$ is the fibre coordinate of the the projection $\pi_1:R\rightarrow M$.
This gives the differential system
\begin{align*}
\tilde D_1&=\{o_1=o_2=o_3=o_4=0\} 
\end{align*}
on $R$.
Let $t$ be the fibre coordinate of the the projection $\pi_2:R\rightarrow J$.
The differential system 
\[
D_2=\{\varpi=\Up=0\}
\]
on $J$ can be pulled back to the differential system $\tilde D_2$ on $R$ by the 1-forms
\begin{align*}
\varpi&=\der y-p\der x-q\der y,\\
\om_2&=\der p+(t^2H_t-2tH+2\int H\der t)\der x+(H-tH_t)\der y,\\
\om_3&=\der q+(H-tH_t)\der x+H_t\der y,\\
\om_4&=\der y-t \der x,
\end{align*} with $\tilde D_2=\{\varpi=\om_2=\om_3=\om_4=0\}$.
There is a coordinate transformation on $R$ that takes $\tilde D_1$ to $\tilde D_2$ and vice versa, giving rise to the double fibration transform. 
Following \cite{Cartan1910}, starting from coordinates $(x,y,z,p,q,t)$ on $R$ associated to the contact distribution $\tilde D_2$, the coordinate transformation to coordinates $(X,Y,Z,P,Q)$ on $M$ is given by
\begin{align*}
X&=-2t,\\
Y&=p+tq+(2\int H \der t-tH)x+Hy,\\
Z&=z-px-qy-(\frac{1}{2}t^2H_t-tH+\int H\der t)x^2-(H-tH_t)xy-\frac{1}{2}H_ty^2,\\
P&=-\frac{1}{2}q-\frac{1}{2}(H-tH_t)x-\frac{1}{2}H_t y,\\
Q&=\frac{1}{4}H_{tt}(y-tx),
\end{align*}
with the fibre coordinate of the projection $\pi_1:R\rightarrow M$ taken to be
\[
L=\frac{1}{8}(H_{ttt}(tx-y)+H_{tt}x).
\]
It follows that
\begin{align*}
\der Y-P\der X&=t\om_2+\om_3,\\
\der P-Q\der X&=-\frac{1}{2}\om_2,\\
\der Z-\frac{4Q^2}{H_{tt}}\der X&=\varpi-y\om_2-x\om_3.
\end{align*}
Upon changing independent variable $t$ to $-\frac{1}{2}X$, we have $H_{XX}=\frac{1}{4}H_{tt}$, so we get the 1-forms encoding the $(2,3,5)$-distribution
\begin{align*}
\der Y-P\der X,\quad 
\der P-Q\der X,\quad 
\der Z-\frac{Q^2}{H_{XX}}\der X,
\end{align*}
where $H(X)$ satisfies Noth's equation (\ref{h-noth}), with the additional 1-form
\[
\der Q-L\der X=\frac{1}{4}H_{tt}\om_4.
\]
The parameter $L$ determines the prolongation of the Lie algebra to $(X,Y,Z,P,Q,L)$.

The converse transformation starting from $\tilde D_1$ on $R$ to coordinates $(x,y,z,p,q,t)$ encoding $\tilde D_2$ is given by taking
\begin{align*}
t&=-\frac{1}{2}X,\\
y&=\frac{1}{H_{XX}}Q-\frac{X}{2}x,\\
q&=-2P+\frac{2H_X}{H_{XX}}Q-xH,\\
p&=Y-PX+(\frac{XH_X-H}{H_{XX}}) Q+\frac{1}{2}(\int (H-XH_X)\der X)x,\\
z&=Z+\frac{H_X}{H_{XX}^2}Q^2-\frac{2}{H_{XX}}PQ+(Y-\frac{H}{H_{XX}}Q)x+\frac{1}{4}(XH+\int (H-XH_X) \der X)x^2.
\end{align*}
We also have
\[
x=\frac{2}{H_{XX}}(H_{XX}L-H_{XXX}Q),
\]
so that
\begin{align*}
\der y-t \der x=\frac{1}{H_{XX}}\der Q-(\frac{x}{2}+\frac{H_{XXX}}{H_{XX}^2}Q)\der X=\frac{1}{H_{XX}}(\der Q-L \der X).
\end{align*}

From the double fibration transform it is in principle possible to obtain the $\frak{g}_2$ Lie algebra of the contact structure independent of the choice of the solution of Noth's equation. The computations are quite complicated and we defer the presentation of the results. To summarize, we have the following.
\begin{proposition}
By the above coordinates given in the double fibration transform, the Lie algebra of vector fields obtained in Theorems \ref{thm1} and \ref{thm2} can be identified with the $\frak{g}_2$ Lie algebra of vector fields obtained in Theorem 4.1 of \cite{r21c} with the parametrization of $H(t)$ given by the solutions in (\ref{noth1}) and (\ref{noth2}) respectively.

\end{proposition}

\end{document}